\documentclass[graybox]{svmult}

\usepackage{type1cm}        
\usepackage{makeidx}         
\usepackage{graphicx}        
\usepackage{multicol}        
\usepackage[bottom]{footmisc}

\usepackage{newtxtext}       %
\usepackage[varvw]{newtxmath}       
\usepackage{mathtools}
\usepackage{enumitem}
\usepackage[utf8]{inputenc}

\newcommand{\bn}{\ensuremath{\mathbf{n}}}
\newcommand{\bv}{\ensuremath{\mathbf{v}}}
\newcommand{\bw}{\ensuremath{\mathbf{w}}}
\newcommand{\bbB}{\ensuremath{\mathbb{B}}}
\newcommand{\bbG}{\ensuremath{\mathbb{G}}}
\newcommand{\bbI}{\ensuremath{\mathbb{I}}}
\newcommand{\bbO}{\ensuremath{\mathbb{O}}}
\newcommand{\bbR}{\ensuremath{\mathbb{R}}}
\newcommand{\bbT}{\ensuremath{\mathbb{T}}}

\newcommand{\dv}[1]{\,{\mathrm d}#1}
\newcommand{\dx}{\dv{x}}
\DeclareMathOperator{\Div}{div}
\DeclareMathOperator{\trace}{Tr}
\newcommand{\abs}[1]{\lvert{#1}\rvert}
\newcommand{\norm}[1]{\|{#1}\|}
\newcommand{\skp}[2]{({#1\, ,\, #2})}
\newcommand{\dualp}[2]{\left<{#1\, ,\, #2}\right>}

\makeindex             

\begin{document}

\title*{On energy-dissipative finite element approximations for rate-type viscoelastic fluids with stress diffusion}
\titlerunning{On energy-dissipative approximations for viscoelastic fluids}
\author{Dennis Trautwein\orcidID{0009-0000-5136-4566}}
\institute{Dennis Trautwein \at University of Regensburg, Universitätsstraße 31, 93053 Regensburg, Germany, \\ \email{dennis.trautwein@ur.de}}
%
%
\maketitle

\abstract*{
We study a fully discrete finite element approximation of a model for unsteady flows of rate-type viscoelastic fluids with stress diffusion in two and three dimensions.
The model consists of the incompressible Navier--Stokes equation for the velocity, coupled with a diffusive variant of a combination of the Oldroyd-B and the Giesekus model for the left Cauchy--Green tensor. 
The discretization of the model is chosen such that an energy inequality is preserved at the fully discrete level. As a consequence, unconditional solvability and stability for the discrete system are guaranteed and the discrete Cauchy--Green tensor is positive definite. Moreover, subsequences of discrete solutions converge to a global-in-time weak solution, as the discretization parameters tend to zero. In the end, we present numerical convergence tests.
}

\abstract{
We study a fully discrete finite element approximation of a model for unsteady flows of rate-type viscoelastic fluids with stress diffusion in two and three dimensions.
The model consists of the incompressible Navier--Stokes equation for the velocity, coupled with a diffusive variant of a combination of the Oldroyd-B and the Giesekus model for the left Cauchy--Green tensor. 
The discretization of the model is chosen such that an energy inequality is preserved at the fully discrete level. Thus, unconditional solvability and stability for the discrete system are guaranteed and the discrete Cauchy--Green tensor is positive definite. Moreover, subsequences of discrete solutions converge to a global-in-time weak solution, as the discretization parameters tend to zero. In the end, we present numerical convergence tests.
}

\section{Introduction}
\label{sec:1}

The work is dedicated to the numerical analysis of a fully discrete finite element approximation of a model for unsteady flows of rate-type viscoelastic fluids with stress diffusion in two and three dimensions, which has been introduced and studied from the analytical point of view in \cite{bathory_2021_viscoelastic_3D}. The model reads as follows.
Let $T>0$ and $\Omega\subset\bbR^d$, $d\in\{2,3\}$, be a bounded Lipschitz domain with boundary $\partial\Omega$. 
For any $(\mathbf{x},t)\in \Omega\times(0,T)$, find the velocity field $\bv(\mathbf{x},t) \in\bbR^d$, the pressure $p(\mathbf{x},t)\in\bbR$ and the left Cauchy--Green tensor $\bbB(\mathbf{x},t)\in\bbR^{d\times d}$, such that
\begin{subequations}
\begin{align}
    \label{eq:v}
    \partial_t \bv + (\bv\cdot\nabla)\bv
    - \eta \Delta \bv
    + \nabla p
    &=  \Div \bbT_e ,
    \\ 
    \label{eq:div}
    \Div \bv &= 0,
    \\
    \label{eq:B}
    \partial_t \bbB + (\bv\cdot\nabla)\bbB
    + \delta_1 (\bbB-\bbI)
    + \delta_2 (\bbB^2-\bbB)
    &= \nabla\bv \bbB 
    + \bbB (\nabla\bv)^\top
    + \lambda \Delta \bbB.
\end{align}
\end{subequations}
Here, the elastic stress tensor is given by 
\begin{align*}
    \bbT_e = 2\mu (1-\beta) (\bbB-\bbI) + 2\mu \beta (\bbB^2-\bbB).
\end{align*}
The shear viscosity and the elastic shear modulus are denoted by $\eta>0$ and $\mu>0$, respectively. The constants $\delta_1, \delta_2 \geq0$ can be interpreted as relaxation parameters and $\beta\in(0,1)$, $\lambda>0$ are given numbers.
We close the system by employing the initial conditions $\bv(0) = \bv_0$ and $\bbB(0) = \bbB_0$ in $\Omega$ and the boundary conditions
\begin{align}
   \label{eq:bc}
    (\bn\cdot\nabla)\bbB &= \bbO, 
    \quad
    \bv= \mathbf{0}
    \quad \text{on } \partial\Omega \times (0,T),
\end{align}
where $\bn$ denotes the outer unit normal to $\partial\Omega$.
Here, $\mathbf{0}\in\bbR^d$ and $\bbO\in\bbR^{d\times d}$ denote the zero vector and zero matrix, respectively, and $\bbI\in\bbR^{d\times d}$ is the unit matrix. 
We note that the evolution equation \eqref{eq:B} preserves the symmetry of $\bbB$ when $\bbB_0$ is symmetric, as \eqref{eq:B} also holds true with $\bbB$ replaced by $\bbB^\top$. 

The viscoelastic system \eqref{eq:v}--\eqref{eq:B} can be seen as a combination of the viscoelastic Oldroyd-B model ($\beta=\delta_2=0$) and the viscoelastic Giesekus model ($\beta=\delta_1=0$). 
In both cases, stress diffusion is included due to $\lambda>0$, which serves as a mathematical regularization. The case $\beta\in(0,1)$ has been proposed in \cite{bathory_2021_viscoelastic_3D} and is essential for an existence result in three dimensions. In particular, for the viscoelastic system \eqref{eq:v}--\eqref{eq:B}, the authors of \cite{bathory_2021_viscoelastic_3D} proved the existence of weak solutions in three dimensions for arbitrarily large time interval and data. 
Here, the Helmholtz free energy 
\begin{align*}
    \psi(\bbB) \coloneqq \mu (1-\beta) ( \trace\bbB - \ln\det\bbB - d ) + \tfrac12 \mu \beta \abs{\bbB-\bbI}^2
\end{align*}
is very helpful in the case $\beta\in(0,1)$, as it provides sufficient regularity for $\bbB$. 
In the case $\beta=0$, additional technical arguments are required to prove the existence of weak solutions, applicable only in two dimensions \cite{barrett_boyaval_2009, barrett_lu_sueli_2017}.

In this work, we study the model \eqref{eq:v}--\eqref{eq:B} from the numerical point of view. 
%
%
The main objective is to consolidate and adapt key findings from the previous works \cite{GKT_2022_viscoelastic, GT_2023_DCDS}, where, for both cases $\beta=0$ and $\beta>0$, the viscoelastic system has been part of more complex models in the context of tumour growth.
%
This approach introduces fresh perspectives since the sole focus lies on the above model, offering a nuanced synthesis of existing results within this reduced context.
In particular, the main focus of this work lies on a systematic strategy for the construction of an energy-dissipative and converging finite element approximation of \eqref{eq:v}--\eqref{eq:B}. 
To enhance comprehension of the mathematical model, we describe the most important properties of the model \eqref{eq:v}--\eqref{eq:B} such as a suitable variational formulation and an energy inequality in Sect.~\ref{sec:2}.
In Sect.~\ref{sec:3}, we present a finite element approximation that is based on the variational formulation and that fulfills a discrete analogue of the energy inequality. We also provide the main arguments how the unconditional existence and stability can be shown. 
Moreover, in Sect.~\ref{sec:4}, we explain how converging subsequences of discrete solutions can be identified such that they converge to a global-in-time weak solution satisfying the variational formulation, as the discretization parameters tend to zero. 
In Sect.~\ref{sec:5}, we present numerical convergence tests.

\section{Variational formulation and energy inequality}
\label{sec:2}
We now state a variational formulation that will be the foundation of our finite element approximation. Note that the authors of \cite{bathory_2021_viscoelastic_3D} have recently proven the existence of weak solutions for a slightly different but equivalent formulation.
In the following, we denote by {\small $\skp{\cdot}{\cdot}_{L^2}$} the $L^2$-inner product and by {\small $\norm{\cdot}_{L^2}$} the $L^2$-norm. 
Moreover, by {\small $\dualp{\cdot}{\cdot}_{H^1}$} and {\small $\dualp{\cdot}{\cdot}_{H^1_{0,\mathrm{div}}}$} we denote the duality pairings between {\small $H^1(\Omega;\bbR^{d\times d}_{\mathrm{S}})$} and {\small $H^1_{0,\mathrm{div}}(\Omega;\bbR^d) = \{\bw\in H^1_0(\Omega;\bbR^d) \mid \Div \bw = 0\}$} with, respectively, their duals, where {\small $\bbR^{d\times d}_{\mathrm{S}}$} denotes the set of symmetric $(d\times d)$-matrices.

\begin{lemma}
Let $(\bv,p,\bbB)$ denote a smooth solution of the system \eqref{eq:v}--\eqref{eq:B} subject to the initial conditions $\bv(0)=\bv_0$, $\bbB(0)=\bbB_0$ with $\bbB_0$ symmetric and the boundary conditions \eqref{eq:bc}. Then, 
\begin{subequations}
\begin{align}
    \label{eq:v_weak} \nonumber
    0 &= \dualp{\partial_t \bv}{\bw}_{H^1_{0,\mathrm{div}}}
    + \tfrac12 \skp{(\bv\cdot\nabla)\bv}{\bw}_{L^2}
    - \tfrac12 \skp{\bv}{(\bv\cdot\nabla)\bw}_{L^2}
    \\
    &\quad
    + \eta \skp{\nabla\bv}{\nabla\bw}_{L^2}
    - \skp{p}{\Div\bw}_{L^2}
    + \skp{\bbT_e}{\nabla\bw}_{L^2},
    \\ \label{eq:div_weak}
    0 &= \skp{\Div \bv}{q}_{L^2},
    \\ \label{eq:B_weak} \nonumber
    0 &= \dualp{\partial_t \bbB}{\bbG}_{H^1}
    - \skp{\bbB}{(\bv\cdot\nabla)\bbG}_{L^2}
    + \delta_1 \skp{\bbB - \bbI}{\bbG}_{L^2}
    \\
    &\quad
    + \delta_2 \skp{\bbB^2 - \bbB}{\bbG}_{L^2}
    - 2 \skp{\nabla\bv}{\bbG\bbB}_{L^2}
    + \lambda \skp{\nabla\bbB}{\nabla\bbG}_{L^2},
\end{align}
\end{subequations}
for a.e.~$t\in(0,T)$ and for all $\bw\in H^1_0(\Omega;\bbR^d)$, $q\in L^2(\Omega)$ and $\bbG\in H^1(\Omega;\bbR^{d\times d}_{\mathrm{S}})$.
\end{lemma}
\begin{proof}
The variational formulation follows by taking the inner product of \eqref{eq:v}--\eqref{eq:B} with the respective test functions and integrating over $\Omega$ and by parts over $\Omega$ in the terms with the Laplacians and in the convective terms.
Moreover, by the symmetry of $\bbB$ and $\bbG$, we note that it holds $\skp{\nabla\bv \bbB + (\nabla\bv)^\top \bbB} {\bbG}_{L^2} = 2 \skp{\nabla \bv}{\bbG \bbB}_{L^2}$, see \cite{barrett_boyaval_2009}. 
\end{proof}

We now present an energy inequality, which results from the particular structure of the variational formulation \eqref{eq:v_weak}--\eqref{eq:B_weak}.
Similarly to the viscoelastic Oldroyd-B model and the viscoelastic Giesekus model, the positive definiteness of $\bbB$ is an essential physical property which is needed to derive an energy inequality \cite{barrett_boyaval_2009, bathory_2021_viscoelastic_3D, GKT_2022_viscoelastic}. 

\begin{lemma}
Let $(\bv, p, \bbB)$ denote a smooth solution of the variational formulation \eqref{eq:v_weak}--\eqref{eq:B_weak}. Then, if $\bbB$ is positive definite, it holds 
\begin{align}
    \label{eq:stability} \nonumber
    &\sup_{t\in(0,T)} \int_\Omega 
    \Big( \tfrac12 \abs{\bv(t)}^2 + \psi(\bbB(t)) \Big) \dx
    + \eta \norm{\nabla\bv}_{L^2(\Omega_T)}^2
    + \mu\lambda \beta \norm{\nabla\bbB}_{L^2(\Omega_T)}^2
    \\
    &\quad \nonumber
    + \mu\lambda \tfrac{1-\beta}{d} \norm{\nabla \ln\det\bbB}_{L^2(\Omega_T)}^2 
    + \mu (1-\beta)\delta_1 \norm{\bbB^{\frac12}-\bbB^{-\frac12}}_{L^2(\Omega_T)}^2
    \\
    &\quad \nonumber
    + \mu 
    \beta\delta_2 \norm{\bbB^{\frac32}-\bbB^{\frac12}}_{L^2(\Omega_T)}^2
    + \mu ((1-\beta)\delta_2 + \beta\delta_1) \norm{\bbB-\bbI}_{L^2(\Omega_T)}^2
    \\
    & \leq \int_\Omega 
    \Big( \tfrac12 \abs{\bv_0}^2 + \psi(\bbB_0) \Big) \dx,
\end{align}
where $\psi(\bbB) \coloneqq \mu (1-\beta) ( \trace\bbB - \ln\det\bbB - d ) + \tfrac12 \mu \beta \abs{\bbB-\bbI}^2$ denotes the Helmholtz free energy of the system and $\Omega_T\coloneqq \Omega\times (0,T)$.
\end{lemma}

\begin{proof}
This follows from choosing the test functions $\bw=\bv$, $q=p$ and $\bbG= \psi'(\bbB) = \mu(1-\beta) (\bbI - \bbB^{-1}) + \mu\beta (\bbB-\bbI)$ in \eqref{eq:v_weak}--\eqref{eq:B_weak} and integrating in time.
We only note
\begin{align}
    \label{eq:kettenregel_formal}
    - \skp{\bbB}{(\bv\cdot\nabla) \psi'(\bbB)}_{L^2}
    &= \skp{\bv}{ - \tfrac12\mu\beta \nabla \abs{\bbB}^2 + \mu (1-\beta) \nabla \ln\det(\bbB^{-1})}_{L^2} = 0,
\end{align}
where the last equality follows from integration by parts and \eqref{eq:div_weak},
see also \cite{barrett_boyaval_2009, GT_2023_DCDS}.
We also note the inequality $-\nabla\bbB : \nabla\bbB^{-1} \geq \tfrac1d \abs{\nabla \ln\det\bbB}^2$, see \cite[Lem.~3.1]{barrett_lu_sueli_2017}. 
\end{proof}

\section{Numerical approximation}
\label{sec:3}

In the remainder of this work, we make the following assumptions. Let $T>0$ and let $\Omega\subset\bbR^d$, $d\in\{2,3\}$, be a convex and 
bounded Lipschitz domain with polygonal (or polyhedral, respectively) boundary $\partial\Omega$. 
We split the time interval $[0,T)$ into equidistant subintervals $[t^{n-1},t^n)$ with $t^n = n \Delta t$ and $t^{N_T}=T$, where $\Delta t = \tfrac{T}{N_T}$, $N_T\in\mathbb{N}$ and $n \in \{0, \ldots, N_T\}$.
We require $\{\mathcal{T}_h\}_{h>0}$ to be a family of conforming partitionings of $\Omega$ into disjoint open simplices $K$ such that $\overline\Omega = \bigcup_{K\in\mathcal{T}_h} \overline K$. 
By $h_K$ we denote the diameter of a simplex $K\in\mathcal{T}_h$ and we set $h=\max_{K\in\mathcal{T}_h} h_K$.
Besides, we assume that $\{\mathcal{T}_h\}_{h>0}$ is shape regular (or non-degenerate). 
Moreover, we assume that the family of meshes $\{\mathcal{T}_h\}_{h>0}$ consists only of non-obtuse simplices, i.e., all dihedral
angles of any simplex in $\mathcal{T}_h$ are less than or equal to $\tfrac\pi2$.

For the discretization of \eqref{eq:v}--\eqref{eq:B}, we introduce the following finite dimensional function spaces: 
\begin{align*}
    \mathcal{S}_h &\coloneqq \{ q_h \in C(\overline\Omega) \mid q_h|_{K} \in \mathcal{P}_1(K;\bbR) \ \forall \,  K\in\mathcal{T}_h \} ,
    \\
    \mathcal{W}_h &\coloneqq \{ \bbB_h \in C(\overline\Omega;\bbR^{d\times d}_{\mathrm{S}}) \mid \bbB_h|_{K} \in \mathcal{P}_1(K; \bbR^{d\times d}_{\mathrm{S}}) \ \forall \,  K\in\mathcal{T}_h \} ,
    \\
    \mathcal{V}_h &\coloneqq \{ \bv_h \in C(\overline\Omega;\bbR^d) \cap H^1_0(\Omega;\bbR^d) \mid \bv_{h}|_{K} \in \mathcal{P}_2(K;\bbR^d) \ \forall \,  K\in\mathcal{T}_h \} .
\end{align*}
Here, $\mathcal{P}_j(K;X)$ denotes the set of polynomials of order $j\in\{1,2\}$ on $K\in\mathcal{T}_h$ with values in $X\in\{\bbR, \bbR^d, \bbR^{d\times d}_{\mathrm{S}} \}$.
Due to the symmetry of $\bbB$, we only consider finite element functions $\bbB_h\in\mathcal{W}_h$ that are symmetric, i.e., $\bbB_h=\bbB_h^\top$. This also reduces the total number of degrees of freedom.
It is well-known that the specific choice for $\mathcal{V}_h \times \mathcal{S}_h$ is inf--sup stable, as the family of meshes $\{\mathcal{T}_h \}_{h>0}$ is shape-regular, see \cite{ern_guermond_2004}. 
Sometimes, an additional but very mild constraint on the mesh is posed due to technical reasons, but it can be dropped if the triangulation is fine enough, see \cite{BBF2013_fem}.  

We recall the standard nodal interpolation operator $\mathcal{I}_h\colon\,C(${\small $ \overline\Omega$}$)\to \mathcal{S}_h$ which is naturally extended to $\mathcal{I}_h\colon\, C(${\small $ \overline\Omega;\bbR^{d\times d}_\mathrm{S}$}$)\to \mathcal{W}_h$, see \cite{barrett_boyaval_2009, GT_2023_DCDS}. 
Moreover, we introduce the following semi-inner products on $C(${\small $ \overline\Omega$}$)$ and $C(${\small $ \overline\Omega;\bbR^{d\times d}_\mathrm{S}$}$)$ by
$\skp{a}{b}_h \coloneqq \int_\Omega \mathcal{I}_h [a b ] \dx$ for any $a,b\in C(${\small $ \overline\Omega$}$;\bbR)$ and $\skp{\mathbb{A}}{\bbB}_h \coloneqq \int_\Omega \mathcal{I}_h [ \mathbb{A}:\bbB ] \dx$ for any $\mathbb{A},\bbB \in C(${\small $ \overline\Omega;\bbR^{d\times d}_\mathrm{S}$}$)$, respectively.
The corresponding semi-norms are both denoted by $\norm{\cdot}_h \coloneqq ${\small $\sqrt{\skp{\cdot}{\cdot}_h}$}.

Now the discrete system of our interest reads as follows.
Suppose that $\bv_h^0 \in \mathcal{V}_{h,\mathrm{div}} \coloneqq \{\mathbf{u}_h\in\mathcal{V}_h \mid \skp{\Div \mathbf{u}_h}{q_h}_{L^2}=0 \, \forall\, q_h\in\mathcal{S}_h\}$ and $\bbB_h^0 \in \mathcal{W}_h$ with $\bbB_h^0$ positive definite are given. Then, for any $n\in\{1, \ldots, N_T\}$, the goal is to find a solution tuple $(\bv_h^n, p_h^n, \bbB_{h}^{n}) \in \mathcal{V}_h \times \mathcal{S}_h  \times \mathcal{W}_h$ with $\bbB_h^n$ positive definite, which satisfies for any test function tuple $(\bw_h, q_h, \bbG_h) \in \mathcal{V}_h \times \mathcal{S}_h \times \mathcal{W}_h$:
\begin{subequations}
\begin{align}
    \nonumber
    \label{eq:v_FE}
    0 &= 
    \tfrac{1}{\Delta t}\skp{\bv_h^n - \bv_h^{n-1}}{\bw_h}_{L^2}
    + \tfrac12 \skp{(\bv_h^{n-1} \cdot \nabla) \bv_h^n}{\bw_h}_{L^2}
    - \tfrac12 \skp{\bv_h^n}{(\bv_h^{n-1} \cdot \nabla) \bw_h}_{L^2}
    \\
    &\quad
    +\eta \skp{\nabla \bv_h^n}{\nabla\bw_h}_{L^2}
    - \skp{p_h^n}{\Div\bw_h}_{L^2}
    + \skp{\mathcal{I}_h \bbT_{e,h}^n}{\nabla\bw_h}_{L^2} ,
    \\
    \label{eq:div_FE}
    0 &= \skp{\Div\bv_h^n}{q_h}_{L^2} ,
    \\
    \label{eq:B_FE}
    \nonumber
    0 &= 
    \tfrac{1}{\Delta t} \skp{\bbB_h^n - \bbB_h^{n-1}}{\bbG_h}_h
    - {\textstyle \sum}_{i,j=1}^d \skp{(\bv_h^{n-1})_i \mathbf{\Lambda}_{i,j}(\bbB_h^n)}{\partial_{x_j} \bbG_h}_{L^2}
    + \delta_1 \skp{\bbB_h^n-\bbI }{\bbG_h}_h 
    \\ 
    &\quad 
    + \delta_2 \skp{(\bbB_h^n)^2 - \bbB_h^n}{\bbG_h}_h 
    - 2 \skp{ \nabla\bv_h^n}{\mathcal{I}_h [ \bbG_h \bbB_h^n ]}_{L^2}
    + \lambda \skp{\nabla\bbB_h^n}{\nabla\bbG_h}_{L^2},
\end{align}
\end{subequations}
where we use the notation $\bbT_{e,h}^n \coloneqq  2 \mu (1-\beta) (\bbB_h^n - \bbI) + 2 \mu \beta \big( (\bbB_h^n)^2 - \bbB_h^n \big)$.

Here, $\{ [\mathbf\Lambda_{i,j}(\bbB_h)]_{k,\ell}\}_{i,j,k,\ell\in\{1,\ldots,d\}}$ is a fourth-order tensor field that has been introduced in \cite{GT_2023_DCDS}. 
The definition of $\mathbf\Lambda_{i,j}(\bbB_h)$ requires $\bbB_h\in\mathcal{W}_h$ to be positive definite and, on any simplex $K\in\mathcal{T}_h$, it holds approximately $\mathbf\Lambda_{i,j}(\bbB_h)|_{K} \approx \delta_{i,j} \, \bbB_h |_{K}$, where $\delta_{i,j}$ denotes the Kronecker symbol.
Moreover, $\mathbf\Lambda_{i,j}(\bbB_h)$ is constructed in a way such that at the fully discrete level we can mimic the identity \eqref{eq:kettenregel_formal},
i.e., for any $\bbB_h\in\mathcal{W}_h$ positive definite and any $\bv_h\in\mathcal{V}_h$ satisfying \eqref{eq:div_FE} it holds
\begin{align}
    \label{eq:Lambda_kettenregel} \nonumber
    &- {\textstyle \sum}_{i,j=1}^d  \skp{(\bv_h)_i \, \mathbf{\Lambda}_{i,j}(\bbB_h)} {\partial_{x_j} \mathcal{I}_h[\beta \bbB_h - (1-\beta) \bbB_h^{-1}] }_{L^2}
    \\
    &= \skp{\bv_h} {\nabla \mathcal{I}_h \big[ -\tfrac12 \beta \abs{\bbB_h}^2 + (1-\beta) \ln\det(\bbB_h^{-1}) \big]}_{L^2} 
    = 0.
\end{align}
For the precise definition of $\mathbf\Lambda_{i,j}(\bbB_h)$, we refer to \cite{GT_2023_DCDS}. We note that ${\mathbf\Lambda}_{i,j}(\bbB_h)$ has originally been introduced in \cite{GT_2023_DCDS} for the case $\beta=\tfrac12$, but the definition can easily be adapted to the case $\beta\in(0,1)$. A similar approach can be found in \cite{barrett_boyaval_2009}. 

We now present the main result for the discrete system \eqref{eq:v_FE}--\eqref{eq:B_FE}.

\begin{theorem}
\label{theorem:existence_FE} 
Let $\eta, \mu, \lambda, \delta_1, \delta_2>0$ and $\beta\in(0,1)$. 
Let $\bv_h^0\in\mathcal{V}_{h,\mathrm{div}}$ and $\bbB_h^0\in\mathcal{W}_h$ be given with $\bbB_h^0$ being positive definite. 
Then, for any $n\in\{1,\ldots,N_T\}$, there exists at least one solution tuple $(\bv_h^n, p_h^n,\bbB_{h}^{n}) \in \mathcal{V}_h \times \mathcal{S}_h \times \mathcal{W}_h$ to the discrete system \eqref{eq:v_FE}--\eqref{eq:B_FE} with $\bbB_{h}^{n}$ being positive definite.
Moreover, all solutions of \eqref{eq:v_FE}--\eqref{eq:B_FE} are stable in the sense that
\begin{align}
    \label{eq:stability_FE} \nonumber
    & \tfrac{1}{2 \Delta t} \norm{\bv_h^n}_{L^2}^2 
    + \tfrac{1}{2 \Delta t} 
    \norm{\bv_h^n - \bv_h^{n-1}}_{L^2}^2
    + \tfrac{1}{\Delta t}  \skp{\psi(\bbB_h^n)}{1}_h
    + \tfrac{\mu}{2 \Delta t}  \norm{\bbB_h^n- \bbB_h^{n-1}}_{h}^2
    \\ \nonumber
    &\quad
    + \eta \norm{\nabla\bv_h^n}_{L^2}^2
    + \mu\beta \Big( 
    \lambda \norm{\nabla \bbB_h^n}_{L^2}^2
    + \delta_1 \norm{\bbB_h^n-\bbI}_h^2
    + \delta_2 \norm{(\bbB_h^n)^{\frac32}-(\bbB_h^n)^{\frac12}}_h^2 \Big)
    \\ \nonumber
    &\quad
    + \mu(1-\beta) \Big(
    \tfrac{\lambda}{d} \norm{\nabla \mathcal{I}_h \ln \det (\bbB_h^n)}_{L^2}^2
    + \delta_1 \norm{(\bbB_h^n)^{\frac12}-(\bbB_h^n)^{-\frac12}}_h^2
    + \delta_2 \norm{\bbB_h^n-\bbI}_h^2
    \Big)
    \\
    &\leq 
    \tfrac{1}{2 \Delta t}  \norm{\bv_h^{n-1}}_{L^2}^2 
    + \tfrac{1}{\Delta t} \skp{\psi(\bbB_h^{n-1})}{1}_h.
\end{align}
\end{theorem}


In the following, we sketch the main arguments for the proof of Theorem \ref{theorem:existence_FE}. The main difficulty is to prove the positive definiteness of $\bbB_h^n$. For that reason, we introduce a regularized system with cut-offs in certain terms where the discrete Cauchy--Green tensor does not necessarily have to be positive definite. 
We also restrict to weakly divergence-free test functions $\bw_h \in \mathcal{V}_{h,\mathrm{div}}$ in the velocity equation. This allows to temporarily forget about the pressure which is reconstructed afterwards using the inf--sup stability of the discrete function spaces for the velocity and the pressure. Note that in practice we would directly compute a solution of \eqref{eq:v_FE}--\eqref{eq:B_FE}.

Let $\delta\in(0,1)$. We define $[s]_\delta\coloneqq \max\{s,\delta\}$ for any $s\in\bbR$ and we extend this definition to matrix-valued functions via eigenvalues, see \cite{barrett_boyaval_2009}.
Suppose that $\bv_{h,\delta}^0 \in \mathcal{V}_{h,\mathrm{div}}$ and $\bbB_{h,\delta}^0 \in \mathcal{W}_h$ with $\bbB_{h,\delta}^0$ not necessarily positive definite are given. Then, for any $n\in\{1,\ldots,N_T\}$, the goal is to find a solution tuple $(\bv_{h,\delta}^n, \bbB_{h,\delta}^{n}) \in \mathcal{V}_{h,\mathrm{div}} \times \mathcal{W}_h$ with $\bbB_{h,\delta}^n$ not necessarily positive definite,
which satisfies for any test function tuple $(\bw_h, \bbG_h) \in \mathcal{V}_{h,\mathrm{div}} \times \mathcal{W}_h$:
\begin{subequations}
\begin{align}
    \nonumber
    \label{eq:v_FE_delta}
    0 &= 
    \tfrac{1}{\Delta t}\skp{\bv_{h,\delta}^n - \bv_{h,\delta}^{n-1}}{\bw_h}_{L^2}
    + \tfrac12 \skp{(\bv_{h,\delta}^{n-1} \cdot \nabla) \bv_{h,\delta}^n}{\bw_h}_{L^2}
    - \tfrac12 \skp{\bv_{h,\delta}^n}{(\bv_{h,\delta}^{n-1} \cdot \nabla) \bw_h}_{L^2}
    \\
    &\quad
    +\eta \skp{\nabla \bv_{h,\delta}^n}{\nabla\bw_h}_{L^2}
    + \skp{\mathcal{I}_h \bbT_{e,h}^{\delta,n}}{\nabla\bw_h}_{L^2} ,
    \\
    \label{eq:B_FE_delta}\nonumber
    0 &= 
    \tfrac{1}{\Delta t} \skp{\bbB_{h,\delta}^n - \bbB_{h,\delta}^{n-1}}{\bbG_h}_h
    - {\textstyle \sum}_{i,j=1}^d \skp{(\bv_{h,\delta}^{n-1})_i \mathbf{\Lambda}_{i,j}^\delta(\bbB_{h,\delta}^n)}{\partial_{x_j} \bbG_h}_{L^2}
    \\ 
    &\quad \nonumber
    + \delta_1 \skp{[\bbB_{h,\delta}^n]_\delta-\bbI }{\bbG_h}_h 
    + \delta_2 \skp{\bbB_{h,\delta}^n [\bbB_{h,\delta}^n]_\delta - [\bbB_{h,\delta}^n]_\delta}{\bbG_h}_h 
    \\
    &\quad
    - \skp{2 \nabla\bv_{h,\delta}^n}{\mathcal{I}_h [ \bbG_h [\bbB_{h,\delta}^n]_\delta ]}_{L^2}
    + \lambda \skp{\nabla\bbB_{h,\delta}^n}{\nabla\bbG_h}_{L^2},
\end{align}
\end{subequations}
where we define $\bbT_{e,h}^{\delta,n} \coloneqq  2 \mu (1-\beta) ([\bbB_h^n]_\delta - \bbI) + 2 \mu \beta \big( \bbB_h^n [\bbB_h^n]_\delta - \bbB_h^n \big)$.
The term ${\mathbf\Lambda}_{i,j}^\delta(\bbB_h)$ is defined similarly to ${\mathbf\Lambda}_{i,j}(\bbB_h)$ and it fulfills an analogue of \eqref{eq:Lambda_kettenregel}, see \cite{GT_2023_DCDS}.

For the regularized system \eqref{eq:v_FE_delta}--\eqref{eq:B_FE_delta}, we have the following result.

\begin{lemma}
\label{lemma:stability_FE_delta} 
Let $\delta\in(0,1)$, $\eta, \mu, \lambda, \delta_1, \delta_2>0$ and $\beta\in(0,1)$. 
Let $\bv_{h,\delta}^0 \in\mathcal{V}_{h,\mathrm{div}}$ and $\bbB_{h,\delta}^0\in\mathcal{W}_h$ be given with $\bbB_{h,\delta}^0$ not necessarily positive definite. 
Then, for any $n\in\{1,\ldots,N_T\}$, there exists at least one solution $(\bv_{h,\delta}^n, \bbB_{h,\delta}^{n}) \in \mathcal{V}_{h,\mathrm{div}} \times \mathcal{W}_h$ to the discrete system \eqref{eq:v_FE_delta}--\eqref{eq:B_FE_delta}. Moreover, all solutions of \eqref{eq:v_FE_delta}--\eqref{eq:B_FE_delta} fulfill
\begin{align}
    \label{eq:stability_FE_delta} \nonumber
    & \tfrac{1}{2 \Delta t} \norm{\bv_{h,\delta}^n}_{L^2}^2 
    + \tfrac{1}{2 \Delta t} 
    \norm{\bv_{h,\delta}^n - \bv_{h,\delta}^{n-1}}_{L^2}^2
    + \tfrac{1}{\Delta t}  \skp{\psi_\delta(\bbB_{h,\delta}^n)}{1}_h
    + \tfrac{\mu}{2 \Delta t}  \norm{\bbB_{h,\delta}^n- \bbB_{h,\delta}^{n-1}}_{h}^2
    \\ \nonumber
    &\quad
    + \eta \norm{\nabla\bv_{h,\delta}^n}_{L^2}^2
    + \mu\beta \Big( 
    \lambda \norm{\nabla \bbB_{h,\delta}^n}_{L^2}^2
    + \delta_1 \norm{ [\bbB_{h,\delta}^n]_\delta -\bbI}_h^2 \Big)
    \\ \nonumber
    &\quad
    + \mu\beta \delta_2 \norm{(\bbB_{h,\delta}^n - \bbI) [\bbB_{h,\delta}^n]_\delta^{\frac12}}_h^2 
    + \mu(1-\beta) 
    \tfrac{\lambda}{d} \norm{\nabla \mathcal{I}_h \ln\det ([\bbB_{h,\delta}^n]_\delta) }_{L^2}^2
    \\ \nonumber
    &\quad
    + \mu(1-\beta) \Big( \delta_1 \norm{[\bbB_{h,\delta}^n]_\delta^{\frac12}-[\bbB_{h,\delta}^n]_\delta^{-\frac12}}_h^2 
    + \delta_2 \norm{[\bbB_{h,\delta}^n]_\delta-\bbI}_h^2 \Big)
    \\
    &\leq 
    \tfrac{1}{2 \Delta t}  \norm{\bv_{h,\delta}^{n-1}}_{L^2}^2 
    + \tfrac{1}{\Delta t} \skp{\psi_\delta(\bbB_{h,\delta}^{n-1})}{1}_h,
\end{align}
where $\psi_\delta(\bbB) \coloneqq \mu (1-\beta) \trace(\bbB - g_\delta(\bbB) - \bbI) + \tfrac\mu2 \beta \abs{\bbB-\bbI}^2$. Here, we define $g_\delta(s) = \ln(s)$ for any $s\geq \delta$ and $g_\delta(s) = \tfrac{s}{\delta} + \ln(\delta) - 1$ for any $s<\delta$, 
and, again, we extend this definition to matrix-valued functions with the help of the eigenvalues, see \cite{barrett_boyaval_2009}.
\end{lemma}
\begin{proof}
The inequality \eqref{eq:stability_FE_delta} follows from direct calculations in \cite[Sect.~2.4.4]{GT_2023_DCDS}. The main idea is to choose $\bw_h = \bv_{h,\delta}^n$ in \eqref{eq:v_FE_delta} and $\bbG_h = \mathcal{I}_h \big[ \mu(1-\beta) (\bbI - [\bbB_{h,\delta}^n]_\delta^{-1}) + \mu\beta (\bbB_{h,\delta}^n - \bbI) \big]$ in \eqref{eq:B_FE_delta}, adding the resulting equations and applying \cite[Lem.~2.5]{GT_2023_DCDS} and \cite[Lem.~2.1]{barrett_boyaval_2009} as well as an analogue of the discrete chain rule \eqref{eq:Lambda_kettenregel} for the regularized case. Moreover, we apply \cite[Lem.~2.9]{GT_2023_DCDS} where the assumption of non-obtuse simplices is used.
The existence of a discrete solution of \eqref{eq:v_FE_delta}--\eqref{eq:B_FE_delta} relies on a combination of the stability estimates \eqref{eq:stability_FE_delta} with a fixed-point argument with, e.g., Brouwer's fixed point theorem. For more details, we refer to, e.g., \cite{barrett_boyaval_2009, GKT_2022_viscoelastic, GT_2023_DCDS}.
\end{proof}

The final step to prove Theorem \ref{theorem:existence_FE} is to pass to the limit $\delta\to0$ in the regularized discrete system \eqref{eq:v_FE_delta}--\eqref{eq:B_FE_delta} and in the \textit{a priori} estimate \eqref{eq:stability_FE_delta}. By choosing $\bv_h^0=\bv_{h,\delta}^0 \in \mathcal{V}_{h,\mathrm{div}}$ and $\bbB_h^0=\bbB_{h,\delta}^0 \in \mathcal{W}_h$, where $\bbB_h^0$ is positive definite, the limit passage $\delta\to0$ can be achieved analogously to \cite[Sect.~2.4.7]{GT_2023_DCDS} and \cite[Thm.~5.2]{barrett_boyaval_2009}.
We note that the positive definiteness of the discrete Cauchy--Green tensor $\bbB_h^n$ from the original system \eqref{eq:v_FE}--\eqref{eq:B_FE} is guaranteed as all non-positive eigenvalues of the discrete Cauchy--Green tensor $\bbB_{h,\delta}^n$ from the $\delta$-regularized system are controlled uniformly in $\delta>0$, which is due to \eqref{eq:stability_FE_delta} and \cite[Lem.~2.5]{GT_2023_DCDS}. 
Moreover, the existence of a pressure $p_h^n \in \mathcal{S}_h$ for the original system \eqref{eq:v_FE}--\eqref{eq:B_FE} follows from the inf--sup stability of $\mathcal{V}_h \times \mathcal{S}_h$.

\section{Convergence of the discrete solutions}
\label{sec:4}
In the following, we present the main steps for the subsequence convergence of solutions of \eqref{eq:v_FE}--\eqref{eq:B_FE} to a weak solution, as $(h,\Delta t) \to (0,0)$.

Let $\bv_0 \in L^2_{\Div}(\Omega;\bbR^d) \coloneqq \{\bw\in L^2(\Omega;\bbR^d) \mid \Div\bw=0 \text{ in } \Omega, \, \bw\cdot\bn=0 \text{ on } \partial\Omega\}$ and $\bbB_0 \in L^2(\Omega;\bbR^{d\times d}_{\mathrm{S}})$ be given with $\bbB_0$ being uniformly positive definite a.e.~in $\Omega$.
We assume that $\bv_h^0 \in \mathcal{V}_{h,\mathrm{div}}$ and $\bbB_h^0 \in \mathcal{W}_h$ are given with $\bbB_h^0$ being positive definite, and that there exists a constant $c_0>0$ that does not depend on $h,\Delta t$, such that
\begin{align}
    \label{eq:stability_initial}
    \tfrac12 \norm{\bv_h^0}_{L^2}^2 
    + \skp{\psi(\bbB_h^0)}{1}_h
    \leq c_0.
\end{align}
In addition, we assume that ($\bv_h^0, \bbB_h^0$) converges weakly to ($\bv_0, \bbB_0$) in $L^2_{\Div}(\Omega;\bbR^d) \times L^2(\Omega;\bbR^{d\times d}_{\mathrm{S}})$.
For specific examples of $\bv_h^0$ and $\bbB_h^0$, we refer to \cite{barrett_boyaval_2009, GKT_2022_viscoelastic, GT_2023_DCDS}.

It follows from \eqref{eq:stability_FE} and \eqref{eq:stability_initial} that there is a $C>0$ not depending on $h,\Delta t$ such that
\begin{align}
\label{eq:stability2_FE}
    & \max_{n\in\{1,\ldots,N_T\}} \Big( \norm{\bv_h^n}_{L^2}^2 
    + \norm{\bbB_h^n}_h^2 \Big)
    + \Delta t {\textstyle \sum}_{n=1}^{N_T} \Big(
    \norm{\nabla\bv_h^n}_{L^2}^2
    + \norm{\nabla\bbB_h^n}_{L^2}^2
    \Big)
    \leq C.
\end{align}
Now, \eqref{eq:stability2_FE} can be used for deriving estimates for the discrete time derivatives of $\bv$ and $\bbB$. More precisely, one can show analogously to \cite[Lem.~3.1]{GT_2023_DCDS} and \cite[Lem.~4.6]{GKT_2022_viscoelastic}, that there exists a constant $C>0$ that is independent of $h,\Delta t>0$, such that
\begin{align}
    \label{eq:higher_order1} 
    & \Delta t {\textstyle \sum}_{n=1}^{N_T} 
    \norm{\tfrac{1}{\Delta t} (\bv_h^n - \bv_h^{n-1})}_{(H^1_{0,\Div})'}^{4/3} 
    +\Delta t {\textstyle \sum}_{n=1}^{N_T} \norm{\tfrac{1}{\Delta t} (\bbB_h^n - \bbB_h^{n-1}) }_{(H^1)'}^{4/3} \leq  C .
\end{align}
The main idea for the proof of \eqref{eq:higher_order1} is to use some specific projections of functions $\bw\in H^1_{0,\Div}(\Omega;\bbR^d)$ and $\bbB\in H^1(\Omega;\bbR^{d\times d}_{\mathrm{S}})$ as test functions in \eqref{eq:v_FE} and \eqref{eq:B_FE}, respectively. 
Then, \eqref{eq:higher_order1} can be derived by using \eqref{eq:stability2_FE} and standard estimates based on Hölder's and Young's inequalities and the Gagliardo--Nirenberg inequality.

Then, we use \eqref{eq:stability2_FE} and \eqref{eq:higher_order1} and standard weak(-$*$) compactness results in order to extract a weakly(-$*$) converging subsequence of discrete solutions, as $(h,\Delta t)\to (0,0)$.
Strongly converging subsequences can be identified by applying the Aubin--Lions theorem 
and a ``compactness by perturbation'' result from \cite{azerad_guillen_2001}, where one needs an additional estimate for the orthogonal Stokes projection of the discrete velocity,
see \cite{GKT_2022_viscoelastic} for more details. 
After that, the limit passage in \eqref{eq:v_FE}--\eqref{eq:B_FE} can be established as in \cite{GT_2023_DCDS}, which is based on the subsequence convergence results together with technical estimates for the numerical errors due to the nodal interpolation operator $\mathcal{I}_h$ and an error estimate for $\mathbf\Lambda_{i,j}(\bbB_h) \approx \delta_{i,j}\, \bbB_h$, $i,j\in\{1,\ldots,d\}$. 
From this, we can deduce the existence of functions $\bv \in L^\infty(0,T;L^2_{\Div}(\Omega;\bbR^d)) \times L^2(0,T;H^1_{0,\Div}(\Omega;\bbR^d))$ and $\bbB \in L^\infty(0,T;L^2(\Omega;\bbR^{d\times d}_{\mathrm{S}})) \times L^2(0,T;H^1(\Omega;\bbR^{d\times d}_{\mathrm{S}}))$ with $\partial_t\bv \in L^{4/3}(0,T; (H^1_{0,\Div}(\Omega;\bbR^d))')$ and $\partial_t\bbB \in L^{4/3}(0,T; (H^1(\Omega;\bbR^{d\times d}_{\mathrm{S}}))')$
that solve the variational formulation \eqref{eq:v_weak}--\eqref{eq:B_weak}. Note that the pressure usually does not enter the variational formulation of the Navier--Stokes equation \eqref{eq:v_weak} at the level of weak solutions, so that we would consider divergence-free functions $\bw\in H^1_{0,\Div}(\Omega;\bbR^d)$ in \eqref{eq:v_weak} instead of $\bw\in H^1_{0}(\Omega;\bbR^d)$.
%
%
One can show analogously to \cite{GT_2023_DCDS} that $\bbB$ is positive definite a.e.~in $\Omega\times (0,T)$ and that the initial conditions $\bv(0)=\bv_0$ and $\bbB(0)=\bbB_0$ are satisfied in the correct sense.
Similarly to \cite{GT_2023_DCDS}, one can recover the energy inequality \eqref{eq:stability} from \eqref{eq:stability_FE}, where the right-hand side of \eqref{eq:stability} is replaced by the constant $c_0$ from \eqref{eq:stability_initial}. 
We note that \eqref{eq:stability} can be correctly recovered from \eqref{eq:stability_FE} if the initial values are given such that $\skp{\psi(\bbB_h^0)}{1}_h \leq \skp{\psi(\bbB_0)}{1}_{L^2}$. 
In general, it is rather hard to construct $\bbB_h^0$ from $\bbB_0$ such that $\skp{\psi(\bbB_h^0)}{1}_h$ is bounded by $\skp{\psi(\bbB_0)}{1}_{L^2}$, which is due to the combination of numerical integration in terms of $\skp{\cdot}{\cdot}_h$ and the logarithmic term in $\psi(\cdot)$.

Finally, this shows that subsequences of discrete solutions of \eqref{eq:v_FE}--\eqref{eq:B_FE} converge to a global-in-time weak solution, as $(h,\Delta t)\to (0,0)$.

\section{Numerical convergence tests}
\label{sec:5}
In this section, we present a convergence study for the numerical scheme \eqref{eq:v_FE}--\eqref{eq:B_FE}, where the implementation is based on the code from \cite{GT_2023_DCDS}. 

We introduce source terms in \eqref{eq:v} and \eqref{eq:B} such that the exact solution $(\tilde\bv, \tilde p, \tilde\bbB)$ of \eqref{eq:v}--\eqref{eq:B} with the corresponding initial and boundary conditions is given by
\begin{align*}
    \tilde\bv(\mathbf{x},t) &= e^{-t}
    \left(\begin{smallmatrix}
        x_1^2 (x_1-1)^2 x_2 (x_2-1) (2x_2-1)\\
        -x_1 (x_1-1) (2x_1-1) x_2^2 (x_2-1)^2
    \end{smallmatrix}\right),
    \quad
    \tilde p(\mathbf{x},t) = e^{-t} (2x_1-1)(2x_2-1),
    \\
    \tilde\bbB(\mathbf{x},t) &= 
    \left(\begin{smallmatrix}
    1 & 0\\ 0 & 1
    \end{smallmatrix}\right)
    + \tfrac{1}{20} e^{-t} \cos(\pi x_1)\cos(\pi x_2)
    \left(\begin{smallmatrix}
       1 & 0\\
        0 & -1
    \end{smallmatrix}\right).
\end{align*}
%
%
%
%
%
%
%
 %
%
%
%
%
We fix the final time $T=0.1$, the domain $\Omega=(0,1)^2$ and the model parameters $\eta = \delta_1 = \mu = \lambda = 1$, $\beta=\tfrac12$ and $\delta_2=0$.
We use the discretization parameters $\Delta t_\ell = \tfrac{1}{5} T \cdot 2^{-\ell}$ and $h_k = 2^{-k}$, where $\ell\in\{3,\ldots,7\}$ and $k\in\{3,\ldots,7\}$, respectively. The uniform spatial meshes are constructed by dividing the domain into $2^k \times 2^k$ squares, and then each square is further subdivided into four identical triangles. 
At each time step, the nonlinear system \eqref{eq:v_FE}--\eqref{eq:B_FE} is solved with a fixed point method until the $\ell^\infty$-norm of the residuum is less than $10^{-12}$. 
Within each fixed point iteration step, the linearized Navier--Stokes subsystem and the linearized Oldroyd-B equation are solved separately with sparse direct solvers. 
On average, it took 2 to 4 iteration steps for the fixed point method to terminate.


For the scheme \eqref{eq:v_FE}--\eqref{eq:B_FE}, we expect first order convergence with respect to time, i.e., $\mathcal{O}(\Delta t)$ for fixed $h$, as the system is discretized with a first order semi-implicit Euler method. In theory, we can only expect first order convergence with respect to space, i.e., $\mathcal{O}(h)$ for fixed $\Delta t$, which is due to an error estimate of order $\mathcal{O}(h)$ for ${\mathbf\Lambda} _{i,j}(\bbB_h) \approx \delta_{i,j} \, \bbB_h$,
for any $i,j\in\{1,\ldots,d\}$ and all $\bbB_h\in\mathcal{W}_h$ with $\bbB_h$ positive definite, see \cite[Lem.~2.8]{GT_2023_DCDS}.

\begin{figure}[t]
\centering
\includegraphics[width=0.95\textwidth,trim={0.2cm 8.5cm 0 8.5cm},clip]{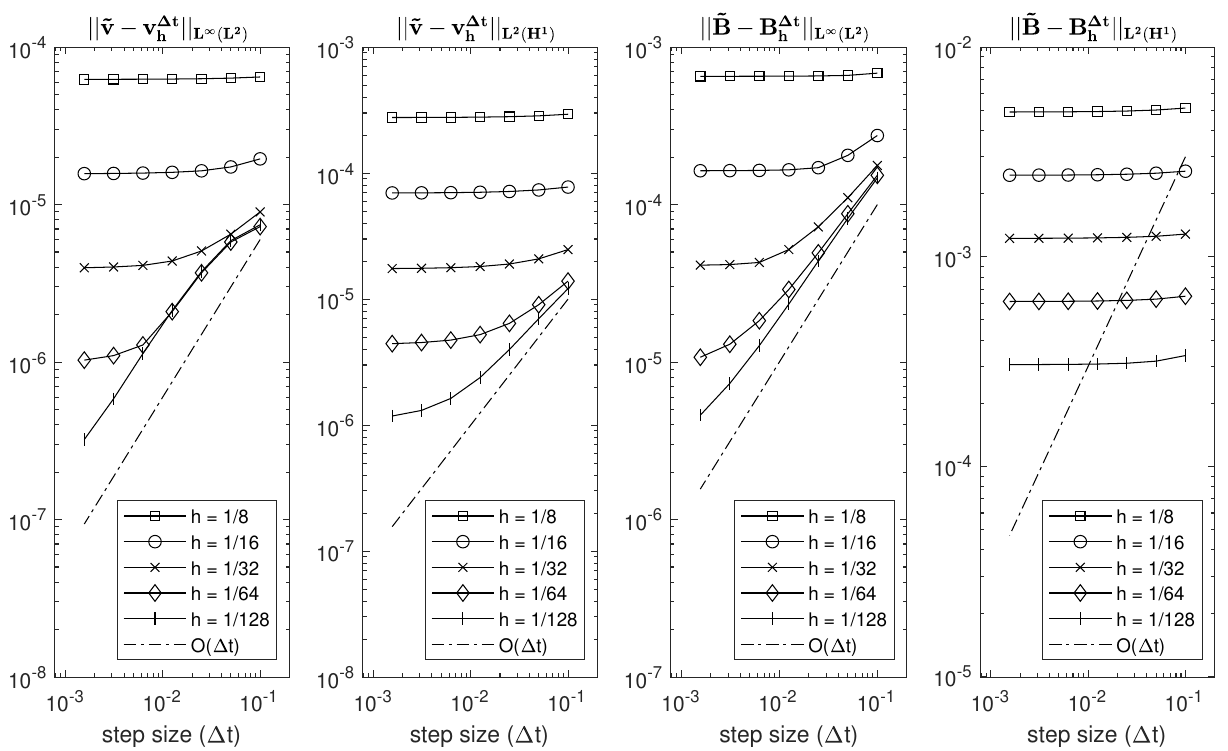}
\caption{Temporal convergence for the velocity and the Cauchy--Green tensor.}
\label{fig:1}      
\end{figure}

\begin{figure}[t]
\centering
\includegraphics[width=0.95\textwidth,trim={0.2cm 8.5cm 0 8.5cm},clip]{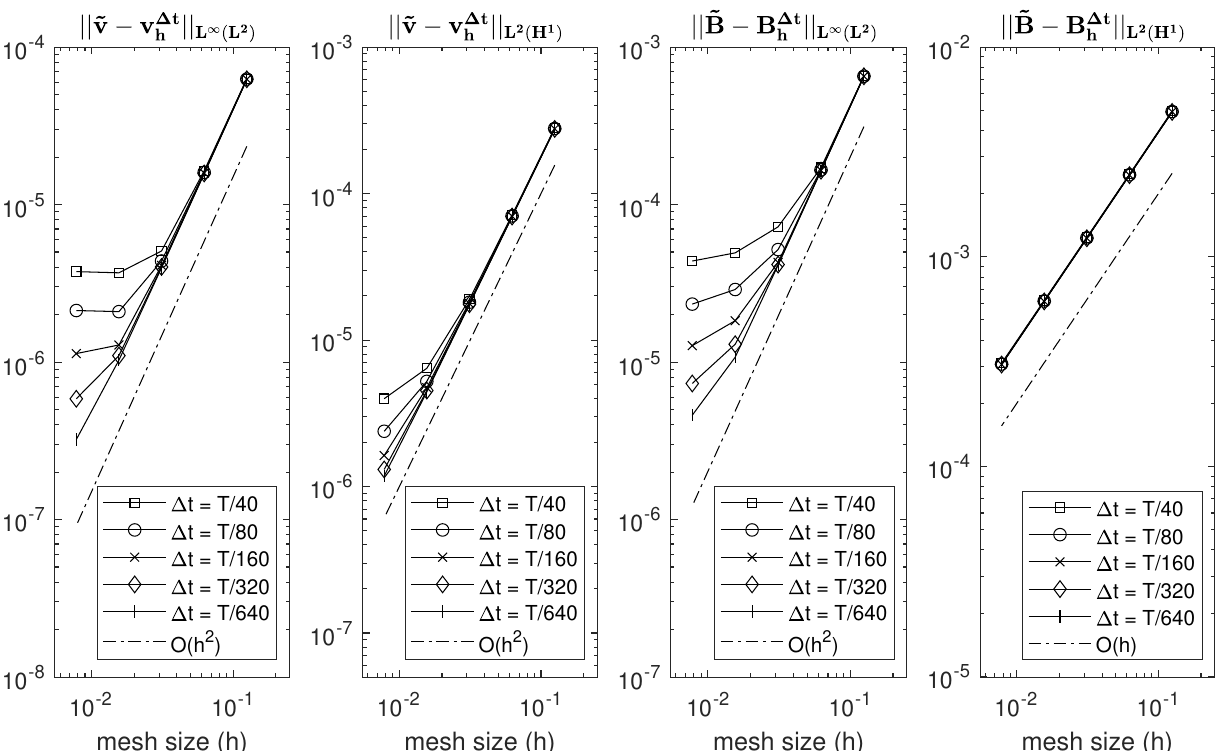}
\caption{Spatial convergence for the velocity and the Cauchy--Green tensor.}
\label{fig:2}       
\end{figure}

In Figs.~\ref{fig:1} and \ref{fig:2}, we visualize the errors between the numerical solution and the exact solution in the $L^\infty(L^2)$ and $L^2(H^1)$ norms. We show the errors for the velocity and the Cauchy--Green tensor against the time step size $\Delta t$ in Fig.~\ref{fig:1}, and against the mesh size $h$ in Fig.~\ref{fig:2}, where lines with distinct symbols represent different mesh and time step sizes, respectively.
In Fig.~\ref{fig:1}, we find a region where the temporal discretisation errors dominate, which matches to the expected $\mathcal{O}(\Delta t)$ convergence, and a region, with small time step size, where the
error curves flatten out, indicating that the spatial discretization errors dominate. In contrast to the expected $\mathcal{O}(h)$ convergence, we observe spatial $\mathcal{O}(h^2)$ convergence for the velocity in the $L^2(H^1)$ norm and for the Cauchy--Green tensor in the $L^\infty(L^2)$ norm, while the $L^2(H^1)$ convergence for the Cauchy--Green tensor is of order $\mathcal{O}(h)$, see Fig.~\ref{fig:2}. These rates are optimal with respect to the approximation properties of the finite element spaces $\mathcal{W}_h$ and $\mathcal{V}_h$. For the $L^\infty(L^2)$ error of the velocity, we also observe $\mathcal{O}(h^2)$ convergence.
These results indicate a rather small influence of the $\mathcal{O}(h)$ approximation of ${\mathbf\Lambda} _{i,j}(\bbB_h) \approx \delta_{i,j} \, \bbB_h$ on the convergence behaviour. However, the presence of ${\mathbf\Lambda} _{i,j}(\bbB_h)$ in the discrete scheme \eqref{eq:v_FE}--\eqref{eq:B_FE} guarantees the unconditional stability, existence and subsequence convergence, as discussed in the previous sections.

\begin{acknowledgement}
The author gratefully acknowledges the support by the Gra\-du\-ier\-ten\-kol\-leg 2339 IntComSin of the Deutsche Forschungsgemeinschaft (DFG, German Research Foundation) -- Project-ID 321821685. 
The author also acknowledges valuable discussions with Harald Garcke.
\end{acknowledgement}

\end{document}